\documentclass[10pt]{article}
\usepackage{amssymb, amsmath, url, graphicx, setspace, geometry}
\usepackage{calrsfs}
\usepackage{wasysym}

\def\3{\subset }
\def\4{\subseteq }
\def\<{\left<}
\def\>{\right>}

\def\bit{\begin{itemize}}
\def\eit{\end{itemize}}
\def\3{\subset }
\def\4{\subseteq }

\def\0{\leqno}

\def\barr{\begin{array}}
\def\earr{\end{array}}

\def\Z{{\rlap{$\kern2pt{\rm Z}$}{\rm Z}\,}}

\title{\bf The second minimum/maximum value of the number of cyclic subgroups of finite $p$-groups}
\author{Mihai-Silviu Lazorec, Rulin Shen and Marius T\u arn\u auceanu}
\date{January 26, 2020}

\begin{document}

\maketitle

\begin{abstract}
Let $C(G)$ be the poset of cyclic subgroups of a finite group $G$ and let $\mathcal{P}$ be the class of $p$-groups of order $p^n$ ($n\geq 3$). Consider the function $\alpha:\mathcal{P}\longrightarrow (0, 1]$ given by $\alpha(G)=\frac{|C(G)|}{|G|}$. In this paper, we determine the second minimum value of $\alpha$, as well as the corresponding minimum points. Further, since the problem of finding the second maximum value of $\alpha$ was completely solved for $p=2$, we focus on the case of odd primes and we outline a result in this regard.       
\end{abstract}

\noindent{\bf MSC (2010):} Primary 20D60; Secondary 20D15, 20D25.

\noindent{\bf Key words:} number of cyclic subgroups, finite $p$-groups, exponent of a group. 

\section{Introduction}

Let $\mathcal{P}$ be the class of $p$-groups of order $p^n$, where $n\geq 3$ is an integer. Let $G\in\mathcal{P}$ and denote by $C(G)$ its poset of cyclic subgroups. Also, for $k\in\lbrace 0,1,\ldots,n\rbrace$,  denote by $c_k(G)$ the number of cyclic subgroups of order $p^k$ of $G$. The literature (see \cite{2, 3}) contains a lot of results concerning the number of subgroups of $G$, most of them having connections with Sylow's theorems which constitute one of the main origins of this research area. Related to the quantities $c_k(G)$, it is worth mentioning that Miller proved in \cite{9} (see also Theorem 1.10 of \cite{2}) that $c_k(G)\equiv 0 \ (mod \ p)$, for all $k\in\lbrace 2,3,\ldots, n\rbrace$, whenever $p$ is an odd prime number. In \cite{10}, the same author showed that the only groups $G\in\mathcal{P}$ satisfying $c_k(G)=p$, for all $k\in\lbrace 2,3,\ldots, n-1\rbrace$, are the abelian group $C_p\times C_{p^{n-1}}$ and the modular $p$-group
$$M(p^n)=\langle x, y \ | \ x^{p^{n-1}}=y^p=1, x^y=x^{1+p^{n-2}}\rangle,$$
the result being valid only for odd primes.

Recently, for a finite group $G$, Garonzi and Lima introduced the ratio $\alpha(G)=\frac{|C(G)|}{|G|}$ in \cite{4}. They completely solved some maxima problems like determining the maximum value of $\alpha(G)$, when $G$ belongs to the class of non-abelian/non-nilpotent/non-supersolvable/non-solvable groups and they also identified the maximum points in each case. We are also interested in this ratio, but we will view it as a function as it follows
$$\alpha:\mathcal{P}\longrightarrow (0,1] \text{ \ given by \ } \alpha(G)=\frac{|C(G)|}{|G|}.$$

In \cite{7} (see Lemma 1), Lazorec and T\u arn\u auceanu proved that
$$\alpha(G)\leq\alpha(C_p^n), \ \forall \ G\in \mathcal{P},$$
and the result can be improved by stating that the equality holds if and only if $exp(G)=p$. In \cite{12}, Richards showed that a finite group of order $n$ has at least $\tau(n)$ cyclic subgroups and the equality holds if and only if the group is isomorphic to $C_n$. A generalization of this result was given by Jafari and Madadi in \cite{5} (see Theorem 2.5). Applying the above result, we have 
$$\alpha(G)\geq \alpha(C_{p^n}), \ \forall \ G\in\mathcal{P} \text{ \ and \ } \alpha(G)=\alpha(C_{p^n})\Longleftrightarrow G\cong C_{p^n}.$$  
In other words, we know the minimum and the maximum values of $\alpha$ and the equality case is also outlined. 

In \cite{1, 8, 11, 15}, the authors deduced the second minimum/maximum value of the number of subgroups on the class $\mathcal{P}$. The aim of this paper is to provide some answers to a similar question that can be formulated as: \textit{"Which is the second minimum/maximum value of $\alpha$ and what can be said about the minimum/maximum points?"} or, equivalently, \textit{"Which is the second minimum/maximum value of the number of cyclic subgroups on the class $\mathcal{P}$ and what can be said about the minimum/maximum points?"}. We investigate this question in the following two sections.

\section{The 2nd minimum value of $\alpha$}

Firstly, we recall the well-known classification of the finite $p$-groups having a cyclic maximal subgroup.\\

\textbf{Theorem 2.1.} \textit{(see Theorem 4.1, \cite{13}, vol. II) Let $G\in\mathcal{P}$ such that $G$ possesses a cyclic maximal subgroup. If $p$ is odd, then $G\cong C_p\times C_{p^{n-1}}$ or $G\cong M(p^n)$, while, if $p=2$, then $G\cong C_2\times C_{2^{n-1}}$ or $G$ is isomorphic to one of the following non-abelian groups:
\begin{itemize}
\item[--] $M(2^n), \ n\geq 4$;
\item[--] the dihedral 2-group $D_{2^n}=\langle x,y \ | \ x^{2^{n-1}}=y^2=1, yxy=x^{-1}\rangle$;
\item[--] the generalized quaternion group $Q_{2^n}=\langle x,y \ | \ x^{2^{n-1}}=y^4=1, yxy^{-1}=x^{2^{n-1}-1}\rangle$;
\item[--] the quasi-dihedral group $QD_{2^n}=\langle x,y \ | \ x^{2^{n-1}}=y^2=1, yxy=x^{2^{n-2}-1}\rangle, \ n\geq 4.$
\end{itemize}}

Using Theorem 2 of \cite{7}, one can compute the number of cyclic subgroups of any finite abelian $p$-group. 
By applying the Inclusion-Exclusion Principle, the number of cyclic subgroups of each of the above non-abelian $p$-groups was determined in \cite{14}. To summarize, we have
\begin{equation}\label{r1}
\begin{cases}|C(C_p\times C_{p^{n-1}})|=|C(M(p^n))|=(n-1)p+2 \\
|C(D_{2^n})|=2^{n-1}+n, \ \ |C(Q_{2^n})|=2^{n-2}+n, \ \ |C(QD_{2^n})|=3\cdot 2^{n-3}+n.
\end{cases}
\end{equation}
Finally, we recall that the number of cyclic subgroups of a finite group $G$ can be computed as follows 
$$|C(G)|=\sum\limits_{x\in G}\frac{1}{\varphi(o(x))},$$ where $\varphi$ is Euler's totient function and $o(x)$ is the order of an element $x$ of $G$. Before stating the main result of this section, we prove the following preliminary fact.\\

\textbf{Lemma 2.2.} \textit{Let $G\in\mathcal{P}$ such that $n\geq 4, G\not\cong C_{p^n}$ and $exp(G)\leq p^{n-2}$. Then
$$|C(G)|>|C(M(p^n))|.$$}

\textbf{Proof.} We proceed by induction on $|G|$. Let $N$ be a maximal subgroup of $G$. Obviously $N\not\cong C_{p^{n-1}}$ since $exp(G)\leq p^{n-2}$. If $exp(N)\leq p^{n-3}$, the inductive hypothesis and (\ref{r1}) imply that
$$|C(G)|=|C(N)|+\sum\limits_{x\in G\setminus N}\frac{1}{\varphi(o(x))}>|C(M(p^{n-1}))|+\sum\limits_{x\in G\setminus N}\frac{1}{\varphi(p^{n-2})}=(n-2)p+2+p^2>|C(M(p^n))|.$$
If $exp(N)=p^{n-2}$, then $N$ has a cyclic maximal subgroup. According to Theorem 2.1 and (\ref{r1}), it follows that $|C(N)|\geq |C(M(p^{n-1}))|$. Hence,
$$|C(G)|=|C(N)|+\sum\limits_{x\in G\setminus N}\frac{1}{\varphi(o(x))}\geq |C(M(p^{n-1}))|+\sum\limits_{x\in G\setminus N}\frac{1}{\varphi(p^{n-2})}\geq (n-2)p+2+p^2>|C(M(p^n))|,$$
which completes the proof.
\hfill\rule{1,5mm}{1,5mm}\\

We are ready to prove the main result of this section which outlines the second minimum value of $\alpha.$\\

\textbf{Theorem 2.3.} \textit{Let $G\in\mathcal{P}$ such that $G\not\cong C_{p^n}$.
\begin{itemize}
\item[i)] If $p$ is odd, then $\alpha(G)\geq\alpha(C_p\times C_{p^{n-1}})=\alpha(M(p^n)).$
\item[ii)] If $p=2$, then 
$\begin{cases} \alpha(G)\geq \alpha(Q_8) &\mbox{, if } n=3 \\ \alpha(G)\geq \alpha(C_2\times C_8)=\alpha(M(16))=\alpha(Q_{16}) &\mbox{, if } n=4 \\ \alpha(G)\geq \alpha(C_2\times C_{2^{n-1}})=\alpha(M(2^n)) &\mbox{, if } n\geq 5\end{cases}.$
\end{itemize}
The equality holds if and only if $G$ is isomorphic to one of the indicated minimum points corresponding to each case.}\\

\textbf{Proof.} Let $G\in \mathcal{P}$. One may use GAP \cite{16} to check the above result for $|G|\in \lbrace 8, 16\rbrace$. Also, for an odd prime $p$ and $n=3$, the classification of finite non-cyclic $p$-groups of order $p^3$ is well-known (see e.g. (4.13), \cite{13}, vol. II). These groups are $C_p\times C_{p^2}$, $C_p^3, M(p^3)$ and $E(p^3)$, where
$$E(p^3)=\langle x,y \ | \ x^p=y^p=[x,y]^p=1, [x,y]\in Z(E(p^3))\rangle.$$
Since $exp(C_p^3)=exp(E(p^3))=p$ and, as we previously mentioned, the number of cyclic subgroups attains its maximum value in this case, our result also holds.      
Consequently, for $p=2$, we directly assume that $n\geq 5$, while, for odd primes, we take $n\geq 4$. Hence, it is sufficient to show that 
$$|C(G)|\geq |C(C_p\times C_{p^{n-1}})|=|C(M(p^n))| \text{ \ and the equality holds iff \ } G\cong C_p\times C_{p^{n-1}} \text{  or  } G\cong M(p^n).$$

We reason by induction on $|G|$. Let $N$ be a maximal subgroup of $G$. If $N\not\cong C_{p^{n-1}}$, by using the inductive hypothesis and (\ref{r1}), we obtain
\begin{align*}
|C(G)|&=|C(N)|+\sum\limits_{x\in G\setminus N}\frac{1}{\varphi (o(x))}\\&\geq |C(C_p\times C_{p^{n-2}})|+\sum\limits_{x\in G\setminus N}\frac{1}{\varphi(p^{n-1})}\\&=|C(C_p\times C_{p^{n-1}})|=|C(M(p^n))|.
\end{align*}
If $N\cong C_{p^{n-1}}$, then $G$ has a cyclic maximal subgroup, so $G$ is isomorphic to one of the groups listed in Theorem 2.1. Then (\ref{r1}) also leads to $|C(G)|\geq |C(C_p\times C_{p^{n-1}})|=|C(M(p^n))|.$

To study the equality case, suppose that $|C(G)|=|C(C_p\times C_{p^{n-1}})|=|C(M(p^n))|$. Then, according to Lemma 2.2, we have $exp(G)=p^{n-1}$, so $G$ is isomorphic to one of the groups listed in Theorem 2.1. Taking into consideration the above equality, it follows that $G\cong C_p\times C_{p^{n-1}}$ or $G\cong M(p^n)$. Since the converse is true, the proof is finished. 
\hfill\rule{1,5mm}{1,5mm}

\section{The 2nd maximum value of $\alpha$}

Let $G\in \mathcal{P}$. We consider the set $\Omega_{\lbrace 1\rbrace}(G)=\lbrace x\in G \ | \ x^p=1\rbrace$ and we denote by $\Omega_1(G)$ the first omega subgroup of $G$, i.e. $\Omega_1(G)=\langle \Omega_{\lbrace 1\rbrace}(G)\rangle$.

If $exp(G)\ne p$ and $p=2$, then Corrolaries 2 and 3 of \cite{4} (see also Lemma 1.4 of \cite{15}), lead to 
$$\alpha(G)\leq \alpha(D_8\times C_2^{n-3}) \text{ \ and \ } \alpha(G)=\alpha(D_8\times C_2^{n-3}) \Longleftrightarrow G\cong D_8\times C_2^{n-3}.$$
Hence, the second maximum value of $\alpha$, as well as the maximum points, are known for $p=2$. In what follows, we are interested in the case of odd primes and we solve the maxima problem while working under the assumption that $\Omega_1(G)\ne G$.\\

\textbf{Theorem 3.1.} \textit{Let $p$ be an odd prime number and let $G\in\mathcal{P}$ such that $exp(G)\ne p$. If $\Omega_1(G)\ne G$, then 
$$\alpha(G)\leq \frac{2p^{n-2}+p^{n-3}+\ldots+p+2}{p^n},$$
and the equality holds if and only if $exp(G)=p^2$ and $\Omega_{\lbrace 1\rbrace}(G)=\Omega_1(G)$ is of index $p$ in $G$.}\\

\textbf{Proof.} Let $exp(G)=p^m$, where $m\geq 2$ is an integer. We have
\begin{align*}
p^n &=1+(p-1)c_1(G)+(p^2-p)c_2(G)+\ldots+(p^m-p^{m-1})c_m(G)\\ &\geq 1+(p-1)c_1(G)+(p^2-p)(c_2(G)+c_3(G)+\ldots +c_m(G)) \\ &=1+(p-1)c_1(G)+(p^2-p)(|C(G)|-1-c_1(G)), 
\end{align*}
so
\begin{equation}\label{r2}
|C(G)|\leq \frac{p^n+p^2-p-1+(p-1)^2c_1(G)}{p^2-p},
\end{equation}
and the equality holds if and only if $exp(G)=p^2$.

Since $\Omega_1(G)\ne G$, it follows that $[G:\Omega_1(G)]\geq p$. Therefore, 
\begin{equation}\label{r3}
|\Omega_{\lbrace 1\rbrace}(G)|\leq |\Omega_1(G)|\leq p^{n-1},
\end{equation}
which leads to
$$c_1(G)\leq\frac{p^{n-1}-1}{p-1}.$$
Using (\ref{r2}) and the above inequality, we obtain
$$|C(G)|\leq\frac{p^n+p^2-p-1+(p-1)(p^{n-1}-1)}{p^2-p}=2p^{n-2}+p^{n-3}+\ldots+p+2,$$
and the conclusion follows. 

Moreover, the equality holds if and only if we have equality in (\ref{r2}) and (\ref{r3}), i.e. $exp(G)=p^2$ and $\Omega_{\lbrace 1\rbrace}(G)=\Omega_1(G)$ is of index $p$ in $G$.
\hfill\rule{1,5mm}{1,5mm}\\ 

Two examples of groups satisfying the equality case of Theorem 3.1 are $M(p^3)\times G_1$, where $|G_1|=p^{n-3}$ and $exp(G_1)=p$, and $C_{p^2}\times G_2$, where $|G_2|=p^{n-2}$ and $exp(G_2)=p.$\\

\textit{Under the same hypotheses as the ones indicated in Theorem 3.1, what can be said if $\Omega_1(G)=G?$} First of all, the result does not hold in this case. For instance, if we take $p=3$ and $n=4$, the wreath product $C_3\wr C_3$ (SmallGroup(81,7)) and the group $E(27)\rtimes C_3$ (SmallGroup(81,9)) satisfy the property $\Omega_1(G)=G$ and they have 29 and 35 cyclic subgroups, respectively. However, the group $M(27)\times C_3$ (SmallGroup(81,13)) has only 23 cyclic subgroups and it checks the equality case of Theorem 3.1. Moreover, note that $\alpha(G)\leq \alpha(E(27)\rtimes C_3)$, 
while the equality holds if and only if $G\cong E(27)\rtimes C_3$.

Secondly, if $\Omega_1(G)=G$, then $G$ can be generated by elements of order $p$. Consequently, $\frac{G}{G'}$ can be generated by elements of order $p$, so $\frac{G}{G'}$ is elementary abelian. Hence $\Phi(G)\subseteq G'$ and, since the converse inclusion is well known, we obtain $G'=\Phi(G)$. If $G'\subseteq Z(G)$, then $G$ is a nilpotent group of class 2 and, since $p$ is odd, it is known that $exp(\Omega_1(G))=p$. This would imply that $exp(G)=p$, a contradiction. Hence, to completely solve the proposed maxima problem a first approach would be to classify the groups $G\in\mathcal{P}$ satisfying $G'=\Phi(G)\not\subseteq Z(G)$ and compare their numbers of cyclic subgroups. A second option would be to find the maximum number of solutions of $x^p=1$ in $G$. This was done in \cite{6} for $p=3$. According to this reference, we have 
$$2c_1(G)\leq 7\cdot 3^{n-2}-1,$$
and the equality holds if and only if $G\cong E(27)\rtimes C_3^2$ (SmallGroup(243,58)) or $G\cong C_3\times E(27)\rtimes C_3$ (SmallGroup(243,53)).
Making the replacements in (\ref{r2}), we obtain
$$|C(G)|\leq \frac{23\cdot 3^{n-3}+1}{2},$$
so 
$$\alpha(G)\leq \frac{23\cdot 3^{n-3}+1}{2\cdot 3^n}.$$ 
The upper bound is attained if and only if one works with a group $G$ isomorphic to $(E(27)\rtimes C_3^2)\times G_1$ or to $ (C_3\times E(27)\rtimes C_3)\times G_1$, where $|G_1|=3^{n-5}$ and $exp(G_1)=3$. \\

We end our paper by suggesting the following open problem:\\

\textbf{Open problem.} \textit{Let $p\geq 5$ be a prime number and let $G\in\mathcal{P}$ such that $\Omega_1(G)=G$. Determine the second maximum value of $\alpha$ in this case.}\\

A solution to this problem along with our results would completely outline the second maximum value of $\alpha$.

\vspace*{3ex}
\small

\begin{minipage}[t]{5cm}
Mihai-Silviu Lazorec \\
Faculty of  Mathematics \\
"Al.I. Cuza" University \\
Ia\c si, Romania \\
e-mail: {\tt mihai.lazorec@student.uaic.ro}
\end{minipage}
\hspace{0.3cm}
\begin{minipage}[t]{5cm}
Rulin Shen \\
Department of  Mathematics \\
Hubei Minzu University\\
Enshi, Hubei, P.R. China \\
e-mail: {\tt rshen@hbmy.edu.cn}
\end{minipage}
\hspace{0.3cm}
\begin{minipage}[t]{5cm}
Marius T\u arn\u auceanu \\
Faculty of  Mathematics \\
"Al.I. Cuza" University \\
Ia\c si, Romania \\
e-mail: {\tt tarnauc@uaic.ro}
\end{minipage}
\end{document}